\documentclass[12pt]{article}
\usepackage{amsfonts}
\usepackage{}
\usepackage{amsthm,amsmath}
\usepackage{amscd}
 \usepackage{amssymb}
 \usepackage{hyperref}
 \usepackage[all]{xypic}
 \usepackage{mathtools}
  \usepackage[numbers,sort&compress]{natbib}
 \usepackage[utf8]{inputenc}
\allowdisplaybreaks[4]
  \newtheorem{lemma}{Lemma}[section]

 \newtheorem{theorem}[lemma]{Theorem}

 \newtheorem{remark}[lemma]{Remark}

\usepackage[utf8]{inputenc} 

\def\includegraphics{}

\usepackage{tikz,xcolor,hyperref}

\definecolor{lime}{HTML}{A6CE39}
\DeclareRobustCommand{\orcidicon}{
\begin{tikzpicture}
\draw[lime, fill=lime] (0,0)
circle[radius=0.16]
node[white]{{\fontfamily{qag}\selectfont \tiny \.{I}D}};
\end{tikzpicture}
\hspace{-2mm}
}
\foreach \x in {A, ..., Z}{%
\expandafter\xdef\csname orcid\x\endcsname{\noexpand\href{https://orcid.org/\csname orcidauthor\x\endcsname}{\noexpand\orcidicon}}
}

\usepackage{amscd}
\usepackage{amsmath,amssymb,mathrsfs}
\usepackage[numbered]{bookmark}
\usepackage{hyperref}
\usepackage[all]{xypic}
\usepackage{times}
\usepackage{inputenc}
\usepackage{titlesec}

\def\inf{\operatorname{inf}}

\allowdisplaybreaks[1]

 \usepackage{amsthm}

\advance\headheight1.15pt
\newenvironment{proof of Theorem 1.1}{{\noindent\it {\bf Proof of Theorem 1.1}}\quad}{\hfill $\square$\par}
\newenvironment{proof of Theorem 1.2}{{\noindent\it {\bf Proof of Theorem 1.2}}\quad}{\hfill $\square$\par}
\newenvironment{proof of Theorem 1.3}{{\noindent\it {\bf Proof of Theorem 1.3}}\quad}{\hfill $\square$\par}
\newenvironment{proof of Theorem 1.4}{{\noindent\it {\bf Proof of Theorem 1.4}}\quad}{\hfill $\square$\par}
\newenvironment{proof of Theorem 1.5}{{\noindent\it {\bf Proof of Theorem 1.5}}\quad}{\hfill $\square$\par}
      \newcommand{\Id}{\int_{\mathbb{D}}}
   \newcommand{\hp}{H^{p}}
\newcommand{\sn}{\sum_{n=0}^{\infty}}  \newcommand{\hu}{\mathcal {H}_{\mu}}
\newcommand{\hd}{H(\mathbb{D})}            
\newcommand{\dd}{\mathbb{D}}                     \newcommand{\pd}{D^{p}_{p-1}}
\newcommand{\B}{\mathcal {B}}
  \newcommand{\h}{\mathcal {H}}

\newcommand{\comment}[1]{}

\begin{document}
\baselineskip=8pt

\title{\fontsize{15}{0}\selectfont  Hilbert matrix operator  on bound analytic functions}

\author{\fontsize{11}{0}\selectfont
  Yuting Guo$^{a}$ and  Pengcheng Tang$^{*,b}$
   \\
   \fontsize{10}{0}\it{
 $^{a}$School of Mathematics and Statistics,
Hunan First Normal University, Changsha, Hunan 410205, China} 
\\ \fontsize{10}{0}\it{  $^{b}$School of Mathematics and Computational Science, Hunan University of Science and Technology,}
 \\ \fontsize{10}{0}\it{ Xiangtan, Hunan 411201, China}
}

\date{}
\maketitle
\thispagestyle{empty}
\begin{center}
\textbf{\underline{ABSTRACT}}
\end{center}
It is well known that the Hilbert matrix  operator $\h$ is bounded from $H^{\infty}$ to the  mean Lipschitz spaces $\Lambda^{p}_{\frac{1}{p}}$ for all $1<p<\infty$.
In this paper, we prove that the range of Hilbert matrix  operator $\h$ acting on   $H^{\infty}$ is contained in certain  Zygmund-type space (denoted by $\Lambda^{1.*}_{1}$),  which is strictly smaller than  $\cap_{p>1}\Lambda^{p}_{\frac{1}{p}}$. We also provide explicit upper and lower bounds for the norm of the Hilbert matrix $\h$ acting from $H^{\infty}$ to  $\Lambda^{1.*}_{1}$. Additionally, we  also characterize the positive Borel measures $\mu$ such that  the generalized Hilbert matrix operator $\hu$ is bounded from $H^{\infty}$ to the Hardy space $H^{q}$.
This part is  a continuation  of the work of  Chatzifountas, Girela and Pel\'{a}ez [J. Math. Anal. Appl. 413 (2014) 154--168] regarding  $\hu$ on  Hardy spaces.

\begin{flushleft}
\end{flushleft}
\begin{flushleft}
\end{flushleft}

\let\thefootnote\relax\footnote{$^*$Corresponding Author}
\let\thefootnote\relax\footnote{ Yuting Guo: rainting2770@163.com }
\let\thefootnote\relax\footnote{ Pengcheng Tang: www.tang-tpc.com@foxmail.com}

\vspace{1cm}

\section{Introduction} \label{Sec:Intro}

 \ \ \  \ \ \ Let $\mathbb{D}=\{z\in \mathbb{C}:\vert z\vert <1\}$ denote the open unit disk of the complex plane $\mathbb{C}$ and $H(\mathbb{D})$ denote the space of all
analytic functions in $\mathbb{D}$.

The Bloch space $\mathcal {B}$ consists of those functions $f\in  H(\mathbb{D}) $ for which
$$
\vert \vert f\vert \vert _{\mathcal {B}}=\vert f(0)\vert +\sup_{z\in \mathbb{D}}(1-\vert z\vert ^{2})\vert f'(z)\vert <\infty.
$$


Let  $0<p\leq\infty$, the classical Hardy space $H^p$  consists of  those functions  $f\in H(\dd)$ for which
$$
||f||_{p}=\sup_{0\leq r<1} M_p(r, f)<\infty,
$$
where
$$
M_p(r, f)= \left(\frac{1}{2\pi}\int_0^{2\pi}|f(re^{i\theta})|^p d\theta \right)^{1/p}, \ 0<p<\infty,
$$
$$M_{\infty}(r, f)=\sup_{|z|=r}|f(z)|.$$

The mixed norm space $H^{p,q,\alpha}$, $0<p,q\leq \infty$, $0<\alpha<\infty$,  is the space of all
functions $f\in \hd$, for which
$$
||f||_{p,q,\alpha}=\left(\int_{0}^{1}M^{q}_{p}(r,f)(1-r)^{q\alpha-1}dr\right)^{\frac{1}{q}}<\infty, \ \mbox{for}\ 0<q<\infty,
$$
and
$$||f||_{p,\infty,\alpha}=\sup_{0\leq r<1}(1-r)^{\alpha}M_{p}(r,f)<\infty.$$

For $t\in \mathbb{R}$,  the fractional derivative of order $t$ of $f\in \hd$ is defined by $D^{t}f(z)=\sum_{n=0}^{\infty}(n+1)^{t}\widehat{f}(n)z^{n}$. If $0<p,q\leq \infty$, $0<\alpha<\infty$, then $H_{t}^{q,p,\alpha}$ is the space of all analytic
functions $f\in \hd$ such that
$$||D^{t}f||_{p,q,\alpha}<\infty.$$

It is a well known fact that if $f\in \hd$, $0<p,q\leq \infty$, $0<\alpha,\beta<\infty$, and $s,t \in \mathbb{R}$ are such that $s-t=\alpha-\beta$, then
$$||D^{s}f||_{p,q,\alpha} \asymp ||D^{t}f||_{p,q,\beta}.$$
Consequently, we get $H_{s}^{p,q,\alpha}\cong H^{p,q,\beta}_{t}$ (see \cite{pm}).

Let  $1\leq p<\infty$ and $0<\alpha\leq 1$, the mean Lipschitz space $\Lambda^p_\alpha$ consists of those functions $f\in H(\dd)$ having a non-tangential limit  almost everywhere such that $\omega_p(t, f)=O(t^\alpha)$ as $t\to 0$. Here $\omega_p(\cdot, f)$ is the integral modulus of continuity of order $p$ of the function $f(e^{i\theta})$. It is  known (see \cite{b1}) that $\Lambda^p_\alpha$ is a subset of $H^p$ and
$$\Lambda^p_\alpha=\left\{f\in H(\dd):M_p(r, f')=O\left(\frac{1}{(1-r^{2})^{1-\alpha}}\right), \ \ \mbox{as}\ r\rightarrow 1\right\}.$$
The space $\Lambda^p_\alpha$ is a Banach space with the norm $||\cdot||_{\Lambda^p_\alpha}$ given by
$$
\|f\|_{\Lambda^p_\alpha}=|f(0)|+\sup_{0\leq r<1}(1-r^{2})^{1-\alpha}M_p(r, f').
$$
It is known (see e.g.  \cite[Theorem 2.5]{BS}) that
$$\Lambda^{p}_{\frac{1}{p}} \subsetneq \Lambda^{q}_{\frac{1}{q}}\subsetneq BMOA \subsetneq \mathcal {B}, \ \ 1<p<q<\infty.$$

For $0<p\leq \infty$,  the Zygmund type space $\mathcal {Z}_{p}$ is the space of  $f\in \hd$ such that
$$||f||_{\mathcal {Z}_{p}}=|f(0)|+|f'(0)|+\sup_{0<r<1}(1-r^{2})M_{p}(r,f'')<\infty.$$

When $p =\infty$, the space $\mathcal {Z}_{\infty}$ is the well known Zygmund space. However,  the space $\mathcal {Z}_{1}$ is closely related to the mean    Lipschitz space $\Lambda^{p}_{\frac{1}{p}}$. For $1<p<\infty$,  with the notations  above,  we see that $\Lambda^{p}_{\frac{1}{p}}\cong H^{p,\infty,1}_{1+\frac{1}{p}}$. On the other hand,  the inclusions between mixed norm spaces(see \cite{ar}) shows that
$$\mathcal {Z}_{1} \cong H_{2}^{1,\infty,1}\cong H_{1+\frac{1}{p}}^{1,\infty,\frac{1}{p}}\subsetneq H^{p,\infty,1}_{1+\frac{1}{p}}\cong\Lambda^{p}_{\frac{1}{p}}.$$
Therefore, the space  $\mathcal {Z}_{1} $  can be regarded as the limit case of $H^{p,\infty,1}_{1+\frac{1}{p}}\cong\Lambda^{p}_{\frac{1}{p}}$ as $p\rightarrow 1$. In view of this point,  we will use the symbol $\Lambda^{1,\ast}_{1}$  instead of  $\mathcal {Z}_{1}$ in the sequel.
Note that
$$\Lambda^{1,\ast}_{1}\subsetneq \Lambda^{p}_{\frac{1}{p}}\subsetneq BMOA \subsetneq\mathcal {B} \ \mbox{ for all}\  1<p<\infty.$$


 Let $\mu$ be a finite positive Borel measure on $[0,1)$ and $n\in \mathbb{N}$. We use $\mu_{n}$ to denote the sequence of order $n$ of $\mu$, that is,  $\mu_{n}=\int_{[0,1)}t^{n}d\mu(t)$. Let $\mathcal {H}_{\mu}$ be the Hankel matrix $(\mu_{n,k})_{n,k\geq 0}$ with entries $\mu_{n,k}=\mu_{n+k}$. The matrix $\mathcal {H}_{\mu}$ induces  an operator on $ H(\mathbb{D}) $ by its action on the Taylor coefficients: $a_{n}\rightarrow\displaystyle{\sum_{k=0}^{\infty}}\mu_{n,k}a_{k},\  n\in \mathbb{N}\cup \{0\}.$
The generalized Hilbert operator $\mathcal {H}_{\mu}$ defined on the spaces $\hd$ of analytic functions in the unit disc $\dd$ as follows:

If $f\in \hd$, $f(z)=\displaystyle{\sum_{n=0}^{\infty}}a_{n}z^{n}$, then
$$
\mathcal {H}_{\mu}(f)(z)=\sum_{n=0}^{\infty}\left(\sum_{k=0}^{\infty}\mu_{n,k}a_{k}\right)z^{n}, \ \ z\in \dd,
$$
whenever the right hand side makes sense and defines an analytic function in $\dd$.
If  $\mu$ is the Lebesgue measure on $[0, 1)$, then  the matrix $\mathcal {H}_{\mu}$ reduces to the classical
Hilbert matrix $\mathcal {H}=(\frac{1}{n+k+1})_{n,k\geq 0}$,  which induces
the classical Hilbert operator  $\mathcal {H}$.

Carleson measures play a key role when we study the generalized Hilbert operators. Recall that if $\mu$ is  a positive Borel measure on $[0, 1)$ and $0 < s < \infty$, then  $\mu$ is
an  $s$-Carleson measure if there exists a positive constant $C$ such that
$$\mu([t,1)) \leq C (1-t)^{s}, \ \ \mbox{for all}\  0\leq t<1 .$$

The study of the Hilbert matrix operator $\mathcal{H}$ on  analytic function spaces was initiated by Diamantopoulos and Siskakis in \cite{dia1}, where they proved that  $\mathcal {H}$ is bounded on Hardy space $H^{p}(1<p<\infty)$, and provided  an upper bound estimate for its norm. Subsequently, Diamantopoulus \cite{dia2} considered the boundedness of $\mathcal {H}$  on the Bergman spaces $A^{p}(2<p<\infty)$ and obtained an upper bound estimate for the norm of $\mathcal {H}$.
 Dostanic, Jevti\'{c} and Vukoti\'{c} extended this work in \cite{dos}, where they provided the exact value of the norm of  $\mathcal {H}$  on the Hardy space $H^{p}(1<p<\infty)$, and determined the precise value of the norm  of  $\mathcal {H}$ on the Bergman space $A^{p}$ for $4<p<\infty$.  However, they left an open problem for the case $2<p<4$ which has been
solved  by Bo\v{z}in and Karapetrovi\'{c} \cite{bok}.
Following these developments, significant research has been devoted to investigating the boundedness of $\mathcal{H}$ and its norm on various analytic function spaces, such as weighted Bergman spaces, mixed norm spaces, Korenblum spaces, and Lipschitz spaces (see \cite{jk,gm,mer,yz,yf,bh,03,ann} and references therein).

In 2012, {\L}anucha, Nowak and Pavlovic \cite{lnp} observed the boundedness of $\mathcal{H}$ from $H^{\infty}$ into $BMOA$.
In fact, it is also true that
$$\mathcal {H}(H^{\infty})\subset \bigcap_{1<p<\infty}\Lambda^{p}_{\frac{1}{p}}\subset BMOA \subset \mathcal {B}.$$

Recently, the  norm of $\h$ from  $H^{\infty}$ to $BMOA$,  to $Q_{p}$ spaces and to the mean Lipchitz spaces $\Lambda^{p}_{\frac{1}{p}}$ has been investigated by Bellavita and  Stylogiannis in \cite{bs}.
 In this note, we will prove that the range of $\mathcal{H}$  acting on $H^{\infty}$ is contained in the Zygmund type space $\Lambda^{1,*}_{1}$. The space $\Lambda^{1,*}_{1}$ is strictly smaller than $\Lambda^{p}_{\frac{1}{p}}$ for any $1<p<\infty$. We also provide  both upper and lower bounds for the norm of  $\mathcal{H}$  from $H^{\infty}$ into $\Lambda^{1,*}_{1}$.

\begin{theorem}
Let $\mu$ be a finite positive  Borel measure on $[0,1)$, then the generalized Hilbert operator  operator $\hu$ is bounded from $H^{\infty}$ to $\Lambda_{1}^{1,\ast}$ if and only if $\mu$ is a Carlenson measure.
\end{theorem}

\begin{theorem}
The norm of $\h$ acting from $H^{\infty}$ to $\Lambda_{1}^{1,\ast}$ satisfies
$$
\frac{3}{2}+\frac{2}{\pi} \leq\| \h \|_{H^{\infty} \rightarrow \Lambda_1^{1,\ast}} \leq \frac{3}{2}+\frac{4}{\pi}.
$$
\end{theorem}

It is clear that the Hilbert matrix operator $\h$ is also bounded from $H^{\infty}$ into $\B$. Here, we shall prove that the exact norm of  $\mathcal{H}$  from  $H^{\infty}$ into $\B$ is equal to 3.
\begin{theorem}
The Hilbert matrix operator  $\h$ is bounded from $H^{\infty}$ to $ \mathcal{B}$ and
$\| \h \|_{H^{\infty}\rightarrow \mathcal{B}}=3.$
\end{theorem}

Widom \cite[Theorem 3.1]{wid} proved that $\mathcal {H}_{\mu}$ is a bounded operator on $H^{2}$ if and only if $\mu$ is a Carleson measure. In 2010, Galanopoulos and Pel\'{a}ez \cite{stu} characterized the positive and finite Borel measures
$\mu$ on $[0, 1)$ for which the generalized Hilbert operator $\mathcal{H}_{\mu}$ is well-defined and bounded on $H^{1}$.  These measures are classified as Carleson-type measures.  In 2014, Chatzifountas, Girela and Pel\'{a}ez  \cite{chat} described the
measures $\mu$ for which $\mathcal {H}_{\mu}$ is a bounded operator from $H^{p}$ to $H^{q}$ for $0<p,q<\infty$. The extreme case $p=q=\infty$  was considered by Girela
and Merch\'{a}n \cite{gire}(see also \cite{bonet}).  However, there are two extreme cases that have not yet been considered: namely, $0<q<p=\infty$ and $0<p<q=\infty$.  Another purpose of this paper is to deal with the extreme case $0< q<p=\infty$. To present our results regarding this question, we will first provide some definitions and notions.

 For $0<p<\infty$, the Dirichlet-type space $\pd$ is the space of  $h\in \hd$ such that
$$||h||_{D^{p}_{p-1}}^{p}=|h(0)|^{p}+ \Id |h'(z)|^{p}(1-|z|)^{p-1}dA(z)<\infty.$$
When $p=2$, the space  $D^{2}_{1}$ is just the  Hardy space $H^{2}$.

 The Hardy-Littlewood space $HL(p)$ consists of those function $f\in \hd$ for which
$$||h||^{p}_{HL(p)}=\sn (n+1)^{p-2}|\widehat{h}(n)|^{p}<\infty.$$

It is well known that
\begin{equation}\label{1}
\pd \subset \hp \subset HL(p),\ \  0<p\leq2,
\end{equation}
\begin{equation}\label{2}
HL(p) \subset \hp \subset \pd,\ \  2\leq p<\infty,
\end{equation}

For $0<q<1$,  let  $B_{q}$ denote the space consisting of those $g\in\hd$ for which
$$||g||_{B_{q}}=\int_{0}^{1}(1-r)^{\frac{1}{q}-2}M_{1}(r,g)dr<\infty.$$
The space $B_{q}$  is consistent with  $H(1,1, \frac{1}{q}-1)$ in the setting of the mixed norm space. The  Hardy space $H^{q}$ is a dense subspace of $B_{q}$ and the two spaces have the same continuous linear functionals \cite{dur}. In \cite{chat}, Chatzifountas, Girela and Pel\'{a}ez  showed that $\mathcal {H}_{\mu}$ is bounded from $H^{p}$ to $B_{q}$ for all $0<p<\infty$ and $0<q<1$, whenever $\mu$ satisfies certain  necessary conditions. Nevertheless, we can know more for $p=\infty$ and $0<q<1$. That is, the operator $\mathcal {H}_{\mu}$ is compact from $H^{p}$ to $B_{q}$ for every  finite positive  Borel measure  $\mu$  on $[0,1)$.

Our main results are stated as follows.
\begin{theorem}
 Let $1\leq q<\infty$ and let $\mu$ be a finite positive  Borel measure on $[0,1)$. Let $Y_{q}\in \{D^{q}_{q-1},H^{q},HL(q)\}$. Then the following statements are equivalent.
 \\ (1) $\mathcal {H}_{\mu}$ is bounded from $H^{\infty}$ to $Y_{q}$.
 \\ (2) $\mathcal {H}_{\mu}$ is compact from $H^{\infty}$ to $Y_{q}$.
 \\ (3) The measure satisfies $\{(n+1)^{1-\frac{2}{q}}\mu_{n}\}_{n=0}^{\infty}\in \ell^{q}$.

\end{theorem}

\begin{theorem}
 Let $0<q<1$ and let $\mu$ be a finite positive  Borel measure on $[0,1)$.  Then $\mathcal {H}_{\mu}$ is compact from $H^{\infty}$ to $B_{q}$.
\end{theorem}

The rest of the paper is organized as follows. Section 2 is devoted to proving Theorems 1.1-1.3, while Section 3 focuses on proving Theorems 1.4 and 1.5.

 Throughout the paper, the letter $C$ will denote an absolute constant whose value depends on the parameters
indicated in the parenthesis, and may change from one occurrence to another. We will use
the notation $``P\lesssim Q"$ if there exists a constant $C=C(\cdot) $ such that $`` P \leq CQ"$, and $`` P \gtrsim Q"$ is
understood in an analogous manner. In particular, if  $``P\lesssim Q"$  and $ ``P \gtrsim Q"$ , then we will write $``P\asymp Q"$.

 \section{The range of  Hilbert operator acting on $H^{\infty}$}

The integral representation of $\hu$ plays a basic role in this work. If $\mu$ is a finite positive  Borel measure on $[0,1)$
 and  $f\in \hd$, we shall write throughout the paper
$$\mathcal {I}_{\mu}(f)(z)=\int_{0}^{1}\frac{f(t)}{(1-tz)}d\mu(t),$$
whenever the right hand side makes sense and defines an analytic function on $\dd$. It
turns out that the operators $\hu$ and $\mathcal {I}_{\mu}$ are closely related.  For instance, if $\mu$ is a Carleson measure, then $\hu (f)=\mathcal {I}_{\mu}(f)$ for all $f\in H^{1}$ \cite{stu}. Since $H^{\infty}\subset H^{1}$, this  is also valid for $f\in H^{\infty}$.

The  following characterization of  Carleson measures on $[0,1)$   is due to  Bao et al. \cite{baoo}.
\begin{lemma}\label{lem2.1}
Suppose  $\beta>0$,  $0\leq  q<s<\infty$ and  $\mu$ is a finite positive  Borel measure on $[0,1)$. Then the following conditions are equivalent:
\begin{enumerate}
  \item [(1)] $\mu$ is a  $s$-Carleson measure;
   \item [(2)]  $$
S_{1}:=\sup_{w\in\dd}\int_{0}^{1}\frac{(1-|w|)^{\beta}}{(1-t)^{q}(1-|w|t)^{s+\beta-q}}d\mu(t)<\infty;
$$
  \item [(3)]$$
S_{2}:=\sup_{w\in\dd}\int_{0}^{1}\frac{(1-|w|)^{\beta}}{(1-t)^{q}|1-wt|^{s+\beta-q}}d\mu(t)<\infty.
$$
\end{enumerate}
\end{lemma}

We  also need the following estimates (see Theorem 1.3 in \cite{lius}).
\begin{lemma}\label{lem2.2}
For $z\in \dd$ and $c\in \mathbb{R}$, define
$$I_{c}(z):=\frac{1}{2\pi}\int^{2\pi}_{0}\frac{1}{|1-ze^{-i\theta}|^{1+c}}d\theta.$$
Then the following statements hold.

(1) If $c<0$, then
$$1\leq I_{c}(z) \leq \frac{\Gamma(-c)}{\Gamma^{2}(\frac{1-c}{2})}.$$

(2) If $c>0$, then
$$1\leq (1-|z|^{2})^{c}I_{c}(z) \leq \frac{\Gamma(c)}{\Gamma^{2}(\frac{1+c}{2})}.$$

(3) If $c=0$, then
$$\frac{1}{\pi}\leq |z|^{2}\left(\log\frac{1}{1-|z|^{2}}\right)^{-1} I_{0}(z) \leq 1.$$

Furthermore, all these inequalities are sharp.

\end{lemma}


\begin{proof of Theorem 1.1}
 If $\mu$ is a Carleson measure, then $\hu(f)=\mathcal {I}_{\mu}(f)$ for all $f\in H^{\infty}$.
By a simple calculation, we have that
\begin{equation}\label{3}
\hu(f)^{\prime\prime}(z)=\int_0^1 \frac{2 f(t) t^2}{(1-t z)^3} d \mu(t). \end{equation}

By (\ref{3}), Fubini's theorem and Lemmas \ref{lem2.1}-\ref{lem2.2}, we obtain
\[ \begin{split}
\sup _{0<r < 1} (1-r^{2}) M_1(r, \hu(f)^{\prime \prime}) &=\sup _{0<r < 1} (1-r^{2}) \frac{1}{2 \pi} \int_0^{2 \pi} \left| \int_0^1 \frac{2 t^2 f(t) d u(t)}{(1-t r e^{i\theta})^3} \right| d \theta \\
& \leqslant\sup _{0<r < 1}(1-r^{2}) \frac{1}{2 \pi} \int_0^{2\pi} \int_0^1 \frac{2 t^2|f (t)|}{|1-t re^ {i\theta}|^3} d \mu(t) d \theta \\
&=  \sup _{0<r < 1}(1-r^{2}) \int_0^1 2t^{2}|f(t)|  \frac{1}{2 \pi} \int_0^{2\pi} \frac{d \theta}{|1-t r e^{i\theta}|^3} d \mu (t) \\
&\leq \frac{2\Gamma(2)}{\Gamma^{2}(\frac{3}{2})}\sup _{0<r < 1}(1-r^{2}) \int_0^1 |f(t)| t^2 \frac{d \mu(t)}{(1-t^{2} r^{2})^2} \\
&\lesssim \|f\|_{H^\infty}\sup _{0<r < 1} \int_0^1 \frac{(1-r^{2})}{(1-t r)^2} d \mu(t) \\
& \lesssim \|f\|_\infty .\\
 \end{split} \]
Therefore, $\mathcal {H}_\mu: H^\infty \rightarrow \Lambda_1^{1,\ast} $ is bounded.\\

On the other hand, if $\hu: H^\infty \rightarrow \Lambda_1^{1,\ast} $ is bounded,
then $ \hu(1)(z)=F_\mu(z)=\sum_{n=1}^{\infty} \mu_n z^{n} \in \Lambda_1^{1,\ast}.$
This implies that
\begin{equation}\label{4}
\sup _{0<r<1}(1-r^{2}) M_1(r, F_\mu^{\prime \prime})<\infty. \end{equation}

By Fej\'{e}r-Riesz inequality and Fubini's theorem, we have
$$
\begin{aligned}
  M_1(r, F_\mu^{\prime \prime})&=\frac{1}{2\pi} \int_0^{2\pi}\left|\int_0^1 \frac{2 t^2 d \mu(t)}{(1-t r e^{i\theta})^2}\right| d \theta \\
& \geq \frac{1}{\pi}\int_0^1 \int_0^1 \frac{2 t^2 d \mu(t)}{(1-t r x)^3} d x \\
&= \frac{1}{\pi}\int_0^1 2 t^2 \int_0^1 \frac{d x}{(1-t r x)^{3}} d \mu(t) \\
& \asymp \int_0^1 \frac{2t^2}{(1-t r)^2} d \mu(t).
\end{aligned}
$$
Using (\ref{4}) and inequalities above, we have that
 \[ \begin{split}
1\gtrsim & \sup _{0<r<1}(1-r^{2}) M_1(r, F_\mu^{\prime \prime})\\
 \gtrsim &\sup_{0<r<1}(1-r^{2}) \int_0^1 \frac{2 t^2}{(1-t r)^2} d\mu(t) \\
 \geq &\sup _{0<r<1}(1-r^{2}) \int_r^1 \frac{2 t^2}{(1-t r)^2} d \mu(t)\\
 \geq &\sup _{\frac{1}{2}<r<1} \frac{(1-r^{2}) 2 r^2}{\left(1-r^2\right)^2} \mu([r, 1)) \\
 \gtrsim &\sup _{\frac{1}{2}<r<1} \frac{\mu([r, 1))}{1-r}. \\
   \end{split} \]
This implies that $\mu$ is a  Carleson measure.
\end{proof of Theorem 1.1}


\begin{proof of Theorem 1.2}
 Let $f(z)=1$, then $\mathcal {H}(f)(0)=\int_{0}^{1}f(t)dt=1$ and $\h(f)'(0)=\int_{0}^{1}tf(t)dt=\frac{1}{2}$.
As shown  previously, we have
 \[ \begin{split}
 M_1\left(r, \h(f)^{\prime \prime}(z)\right)  & \geqslant \frac{1}{\pi} \int_0^1 \int_0^1 \frac{2 t^2}{(1-trx)^3} d t d x \\
& =\frac{1}{\pi} \int_0^1 \int_0^1 2 \sum_{n=0}^{\infty} \frac{\Gamma(3+n)}{\Gamma(n+1)\Gamma(3)} t^{n+2} r^n x^n  d t d x \\
& =\frac{1}{\pi} \sum_{n=0}^{\infty} \frac{n+2}{n+3} r^n.
 \end{split} \]
For $0<r<1$, is is easy to compute that
 \[ \begin{split}
 \sum_{n=0}^{\infty} \frac{n+2}{n+3} r^n &=   \sum_{n=0}^{\infty} r^n-\frac{1}{r^{3}}\sum_{n=0}^{\infty}\frac{r^{n+3}}{n+3}\\
  &=\frac{1}{1-r}-\frac{1}{r^3}\left(\log \frac{1}{1-r}-r-\frac{r^2}{2}\right).
 \end{split} \]
This yields  that
$$
\begin{aligned}
 \|\h(f)\|_{\Lambda_1^{1,\ast}}&=|\h(f)(0)|+|(\h (f))^{\prime}(0)|+\sup _{0<r<1}(1-r^{2}) M_1(r, \h(f)^{\prime \prime}) \\
& =\frac{3}{2}+\sup_{0<r<1}(1-r^{2}) M_1(r, \h(f)^{\prime \prime}) \\
&\geq \frac{3}{2}+\frac{1}{\pi} \sup _{0<r<1}(1-r^{2})\left[\frac{1}{1-r}-\frac{1}{r^3}\left(\log \frac{1}{1-r}-r-\frac{r^2}{2}\right)\right] \\
&=\frac{3}{2}+\frac{1}{\pi} \sup_{0<r<1}(1+r)\left[1-\frac{(1-r)}{r^3}\left( \log\frac{1}{1-r}-r-\frac{r^2}{2}\right)\right]. \\
\end{aligned}
$$
Let $$F(r)=1-\frac{(1-r)}{r^3}\left( \log\frac{1}{1-r}-r-\frac{r^2}{2}\right), \ 0< r<1.$$
After careful calculations, we obtain
$$F'(r)=\frac{\frac{r^{2}}{2}-3r-2r\log\frac{1}{1-r}+3\log\frac{1}{1-r}}{r^{4}}.$$
To show that $F(r)$ is increasing on the interval  $(0,1)$, it is suffices to prove that  $\psi(r)=\frac{r^{2}}{2}-3r-2r\log\frac{1}{1-r}+3\log\frac{1}{1-r}>0$ on  $(0,1)$.
Now, it is easy to check that
$$\psi'(r)=r-1+\frac{1}{1-r}-2\log\frac{1}{1-r} \ \mbox{and}\ \psi''(r)=\frac{r^{2}}{(1-r)^{2}}.$$
Since $\psi(0)=\psi'(0)=0$ and $\psi''(r)>0$ for all $0<r<1$, this means that $\psi(r)>0$ for all $r\in (0,1)$. So we conclude that $F(r)$ is  monotonically increasing on the interval  $(0,1)$. This also implies that $(1+r)\left[1-\frac{(1-r)}{r^3}\left( \log\frac{1}{1-r}-r-\frac{r^2}{2}\right)\right]$ is increasing on  $(0,1)$.

By L'H\"{o}pital's rule we have that
$$
\sup _{0<r <1} (1+r)F(r) =2\lim _{r \rightarrow 1^{-}} \left[1-\frac{1-r}{r^3}\left(\log \frac{1}{1-r}-r-\frac{r^2}{2}\right)\right]=2.
$$
Therefore,  we get $\|\h(f)\|_{\Lambda_1^{1,\star}}\geq \frac{3}{2}+\frac{2}{\pi}$.

On the other hand, for  any $f\in H^{\infty}$, we have
$$|\h (f)(0)|=\left|\int_0^1 f(t) d t\right|\leqslant\|f||_{H^{\infty}} \int_0^1 d t=\| f \|_{ \infty},$$
and
$$
|\h(f)^{\prime}(0)|=\left|\int_0^1 t f(t) d t\right|
\leqslant\|f\|_ \infty \int_0^1 t d t=\frac{1 }{2}\|f\|_ \infty .
$$
By the definition of $\Lambda^{1,*}_{1}$, we get
$$
\begin{aligned}
||\h(f)||_{\Lambda_1^{1,\ast}}&=|\h f(0)|+ |\h(f)^{\prime}(0)|+\sup _{0<r<1}(1-r^{2}) M_1(r, \h(f)^{\prime \prime}) \\
& \leqslant \frac{3}{2} ||f||_{\infty}+\sup _{0<r<1}(1-r^{2}) M_1\left(r, \h(f)^{\prime \prime}\right). \\
\end{aligned}
$$
As the proof of Theorem 1.1 shows, we have
 \[ \begin{split}
 M_1\left(r, \h(f)^{\prime \prime}\right)
 & \leq  ||f||_{H^{\infty}} \int_0^1 2 t^2 \frac{\Gamma(2)}{\Gamma^{2}(\frac{3}{2})}\frac{1}{(1-t^2 r^2)^2}dt\\
 &=||f||_{H^{\infty}}\frac{1}{\pi} \int_0^1 \frac{8t^{2}}{(1-t^2 r^2)^2}dt\\
 &=|| f||_ \infty \frac{1}{\pi} \int_0^1 \frac{4 \rho^{\frac{1}{2}}}{(1-\rho r^2)^{2}} d \rho.
 \end{split} \]
Using above inequalities, we obtain that
\[ \begin{split}
&\sup _{0<r<1}(1-r^{2}) M_1\left(r, \h(f)^{\prime \prime}\right) \\
\leq &||f||_{H^{\infty}}\frac{1}{\pi}\sup _{0<r<1}(1-r^{2}) \int_0^1 \frac{4 \rho^{\frac{1}{2}}}{\left(1-\rho r^2\right)^2} d \rho   \\
 =&||f||_{H^{\infty}} \frac{4}{\pi} \sup_{0<r<1} (1-r^{2}) \int_0^1 \sum_{n=0}^{\infty}(n+1) \rho^{n+\frac{1}{2}} r^{2 n} d \rho\\
 =&||f||_{H^{\infty}} \frac{4}{\pi}\sup _{0<r<1}(1-r^{2}) \sum_{n=0}^{\infty} \frac{n+1}{n+3/2} r^{2n} \\
=&||f||_{H^{\infty}} \frac{4}{\pi}\sup _{0<r<1}(1-r^{2}) \left[\sum_{n=0}^{\infty}  r^{2n}- \sum_{n=0}^{\infty}\frac{r^{2n}}{2n+3} \right]\\
=&||f||_{H^{\infty}} \frac{4}{\pi}\sup _{0<r<1}(1-r^{2}) \left[\frac{1}{1-r^{2}}- \frac{\frac{1}{2}\log\frac{1+r}{1-r}-r}{r^{3}}\right]\\
= & \frac{4}{\pi}||f||_{H^{\infty}}  \left(1-\inf_{0<r<1}\frac{(1-r^{2})(\frac{1}{2}\log\frac{1+r}{1-r}-r)}{r^{3}}\right)\\
=&\frac{4}{\pi} ||f||_{H^{\infty}} .\\
 \end{split} \]
Therefore, $\|\h\|_{H^{\infty} \rightarrow \Lambda_1^{1,\ast}} \leq \frac{3}{2}+\frac{4}{\pi}.$
\end{proof of Theorem 1.2}
\begin{remark}
In \cite{ant}, the author prove that
$$C=\sup_{z\in\dd}\int_{\dd}\frac{2(1-|z|^{2})|w|}{|1-z\overline{w}|^{3}}dA(w)=\sup _{0<r<1}\int_{\dd}\frac{2(1-r^{2})|w|}{|1-r\overline{w}|^{3}}dA(w)=\frac{8}{\pi}.$$
Using this result, we may  easily obtain  an upper bound estimate for the norm of $\h$ from $H^{\infty}$ to $\Lambda_1^{1,\ast}$. As above shows,
 \[ \begin{split}
 M_1\left(r, \h(f)^{\prime \prime}\right)
 & \leq ||f||_{H^{\infty}}\int_0^1 2 t^2 \left(\frac{1}{2 \pi} \int_0^{2 \pi} \frac{d \theta}{|1-t r e^{i\theta}|^3}\right) dt\\
&=||f||_{H^{\infty}} \int_{\dd}\frac{2|w|}{|1-r\overline{w}|^{3}}dA(w).
 \end{split} \]
 This implies that
$$\sup _{0<r<1}(1-r^{2}) M_1\left(r, \h(f)^{\prime \prime}\right) \leq ||f||_{H^{\infty}} \sup _{0<r<1}\int_{\dd}\frac{2(1-r^{2})|w|}{|1-r\overline{w}|^{3}}dA(w)
 = \frac{8}{\pi} ||f||_{H^{\infty}} .$$
\end{remark}

The Hilbert matrix operator $\h$ is bounded from $H^{\infty}$ to $BMOA$. It is well known that $BMOA\subsetneq \B$, and hence the Hilbert matrix operator $\h$ is also bounded from $H^{\infty}$ to $\B$. Here, we give the exact norm of  Hilbert matrix operator acting from  $H^{\infty}$ into $\B$.

\begin{proof of Theorem 1.3}
 Let $ f \in H^{\infty} $ and $f \not \equiv 0$. Then, using  the integral form of $\h(f)$ and after a simple calculation, we obtain
$$
\begin{aligned}
\| \h (f)\|_\mathcal{B}&=|H f(0)|+\sup _{z \in \mathbb{D}}\left(1-|z|^2\right)\left|\h^{\prime} (f)(z)\right| \\
& =\left|\int_0^1 f(t) d t\right|+\sup _{z \in \mathbb{D}}\left(1-|z|^2\right)\left|\int_0^1 \frac{f(t) t}{(1-t z)^2} d t\right| \\
& \leq\|f\|_{H^{\infty}}\left(1+\sup _{z \in \mathbb{D}}\left(1-|z|^2\right) \int_0^1 \frac{t}{(1-t|z|)^2} d t \right). \\
\end{aligned}
$$
On the other hand,   we choose the test function $h(z)=1$. Then, $h\in H^{\infty}$  and $||h||_{\infty}=1$.
 Thus,
  \[ \begin{split}
 \| \h \|_\mathcal{B}&\geq ||\h(h)||_\mathcal{B}=1+\sup_{z \in \mathbb{D}}\left(1-|z|^2\right)\left|\int_0^1 \frac{ t}{(1-t z)^2} d t\right| \\
  & \geq 1+\sup_{0\leq x<1}(1-x^{2})\int_{0}^{1}\frac{t}{(1-tx)^{2}}dt.
    \end{split} \]
Therefore, $$\| \h \|_{H^{\infty}\rightarrow\mathcal{B}}=1+\sup _{0\leq x<1}\left(1-x^2\right) \int_0^1 \frac{t}{(1-tx)^2} d t.$$
By making a change of variables  $t=\frac{1-s}{1-xs}$, we have

$$
\begin{aligned}
 &\sup_{0\leq x<1}(1-x^2) \int_0^1 \frac{t}{(1-tx)^2} d t \\
 =&\sup_{0\leq x<1} \left(1-x^2\right) \int_0^1 \frac{1-s}{1-xs} \left(\frac{1-x s}{1-x}\right)^2  \frac{1-x}{(1-x s)^2} d s \\
 =&\sup_{0\leq x<1}(1+x) \int_0^1 \frac{1-s}{1-xs} d s. \\
\end{aligned}
$$
Let $$ G(x)=(1+x) \int_0^1 \frac{1-s}{1-xs} d s, \ \ x\in [0,1),$$
then
\[ \begin{split}
G(x)&=(1+x)\int_0^1 \sum_{n=0}^{\infty} x^n s^n(1-s)ds\\
&=(1+x) \sum_{n=0}^{\infty} \int_0^1 s^n(1-s)ds\  x^n\\
&=(1+x) \sum_{n=0}^{\infty} B(n+1,2)\ x^n\\
&=(1+x)\left(\sum_{n=0}^{\infty}\frac{x^n}{n+1}-\sum_{n=0}^{\infty} \frac{x^n}{n+2}\right)\\
&=(1+x)\left[\frac{1}{x} \log \frac{1}{1-x}-\frac{1}{x^2}\left(\log \frac{1}{1-x}-x\right)\right].
 \end{split} \]
For fixed $s$ in $[0, 1]$, $(1-xs)^{-1}$ is monotonically increasing with respect to $x$ in $[0, 1)$. It follows that $G(x)$ is monotonically increasing in $[0, 1)$ and hence
$$\sup_{0\leq x<1}G(x)=\lim_{r\rightarrow 1^{-}}G(x)=2.$$
So we have that
 $$\| \h  \|_{H^{\infty}\rightarrow \mathcal{B}}=3.$$
\end{proof of Theorem 1.3}

\begin{remark}
The  upper bound for the norm of  $\h$ from $H^{\infty}$ to $\B$ can also be proved in the following way, which is taken from \cite{bs}. Note that
$$\mathcal {H}(f)'(z)=\int_{0}^{1}\frac{tf(t)}{(1-tz)^{2}}dt.$$
The convergence of the integral and the analyticity of the function $f$ guarantee that we can change the path of integration to
$$\gamma(t)=\frac{t(1-z)}{1-tz},\ \ 0\leq t\leq 1.$$
Therefore, we have that
$$\mathcal {H}(f)'(z)=\frac{1}{1-z}\int_{0}^{1}\frac{t}{1-(1-t)z}f(\frac{t}{1-(1-t)z})dt.$$
Since $\psi_{t}(z)=\frac{t}{1-(1-t)z}$ maps the unit disc into itself for each $0\leq t<1$, it follows that $$|\psi_{t}(z)f(\psi_{t}(z))|\leq ||f||_{\infty}.$$
So we can  rewrite $\mathcal {H}(f)'$ as
\begin{equation}\label{5}
\mathcal {H}(f)'(z)=\frac{g(z)}{1-z},
\end{equation}
where $g\in H^{\infty}$ and $||g||_{\infty}\leq ||f||_{\infty}$.

Now, using (\ref{5}) we get
\[ \begin{split}
\| \h (f)\|_\mathcal{B} &=|\h f(0)|+\sup _{z \in \mathbb{D}}\left(1-|z|^2\right)\left| \h^{\prime} (f)(z)\right| \\
& \leq ||f||_{\infty}+ ||g||_{\infty}\sup _{z \in \mathbb{D}}\frac{(1-|z|^{2})}{|1-z|}\\
& \leq ||f||_{\infty}+ ||f||_{\infty}\sup _{z \in \mathbb{D}}\frac{(1-|z|^{2})}{(1-|z|)}\\
&\leq 3 ||f||_{\infty}.
 \end{split} \]

\end{remark}

Recall that the Ces\`{a}ro operator $\mathcal {C}$  is defined in $\hd$ as follows: If $f\in \hd$, $f(z)=\sum_{n=0}^\infty a_nz^n$, then
 $$
\mathcal {C}(f)(z)=\sum_{n=0}^\infty\left(\frac{1}{n+1}\sum_{k=0}^n a_k\right)z^n=\int_{0}^{1}\frac{f(tz)}{1-tz}dt, \  z\in\dd.
$$
The study of the Ces\`{a}ro operator $\mathcal {C}$  on  various spaces of analytic functions has a fairly long time and the integral form  of  $\mathcal {C}$  is closely related to the Hilbert operator $\h$.
 In \cite{ces5}, Danikas and Siskakis proved that $\mathcal{C}$ is bounded from   $H^{\infty}$ to $BMOA$, and
$\|\mathcal{C}\|_{H^{\infty}\rightarrow BMOA}=1+\frac{\pi}{\sqrt{2}}$. Since $BMOA\subsetneq \mathcal{B}$, the  Ces\`{a}ro operator $\mathcal{C}$ is also bounded from   $H^{\infty}$ to $\mathcal{B}$. Following the above arguments, it is easy to obtain that $$\|\mathcal{C}\|_{H^{\infty}\rightarrow \B}=3.$$

 \section{  Hilbert operator acting from  $H^{\infty}$ to Hardy spaces}



We begin with some preliminary results that will be used repeatedly throughout the  rest of the paper. The first  lemma   provides  a characterization of $L^{p}$-integrability of power series
with nonnegative coefficients.  For a proof, see \cite[Theorem 1]{1983}.

\begin{lemma}\label{lem3.1}
Let  $0<\beta,p<\infty$, $\{\lambda_{n}\}_{n=0}^{\infty}$  be a sequence of  non-negative  numbers. Then
$$\int_{0}^{1}(1-r)^{p\beta-1}\left(\sn \lambda_{n}r^{n}\right)^{p}dr\asymp \sn 2^{-np\beta}\left(\sum_{k\in I_{n}}\lambda_{k}\right)^{p}$$
where  $I_{0}=\{0\}$, $I_{n}=[2^{n-1},2^{n})\cap \mathbb{N}$ for $n\in \mathbb{N}$.
\end{lemma}

The following result can be found in \cite{I1} and hence its proofs is omitted.
\begin{lemma}\label{lem3.2}
Let $1<q <\infty$ and $Y_{q}\in \{D^{q}_{q-1},H^{q}, HL(q)\}$. Suppose $f(z)=\sum_{n=0}^{\infty}a_{n}z^{n}\in \hd$ and the sequence $\{a_{n}\}_{n=0}^{\infty}$  is non-negative decreasing, then
$f\in Y_{q}$ if and only if $$\{(n+1)^{1-\frac{2}{q}}a_{n}\}_{n=0}^{\infty}\in \ell^{q}.$$
\end{lemma}

\begin{lemma}\label{lem3.3}
 Let $\mu$ be a finite positive Borel measure on $[0,1)$. Let $\{f_{k}\}_{k=1}^{\infty}\subset H^{\infty}$ such that  $\sup_{k\geq 1}||f_{k}||_{\infty}<\infty$ and  $f_{k}\rightarrow 0$ uniformly on compact subsets
of $\dd$. Then $$\lim_{k\rightarrow \infty}\int_{0}^{1}|f_{k}(t)|d\mu(t)=0.$$
\end{lemma}
\begin{proof}
Since $\mu$ is a finite positive Borel measure on $[0,1)$, for any $\varepsilon>0$,  there exists a  $\delta_{0}\in (0,1)$ such that $\mu([\delta_{0},1))<\varepsilon$.
By the hypothesis there exists  $k_{0}\in \mathbb{N}$ such that
$$|f_{k}(t)|<\varepsilon,\ \ \mbox{if}\ k\geq k_{0}\ \mbox{and}\ 0\leq t\leq \delta_{0}.$$
We see that for $k\geq k_{0}$, we have
 \[ \begin{split}
 \int_{0}^{1}|f_{k}(t)|d\mu(t)&=\int_{0}^{\delta_{0}}|f_{k}(t)|d\mu(t)+\int_{\delta_{0}}^{1}|f_{k}(t)|d\mu(t)\\
&\lesssim \varepsilon+\sup_{k\geq 1}||f_{k}||_{\infty}\mu([\delta_{0},1))\\
&\lesssim \varepsilon.
   \end{split} \]
   The proof is complete.
   \end{proof}

\begin{proof of Theorem 1.4} It is suffices to prove that $(1)\Rightarrow (3)$ and $(3)\Rightarrow (2)$.

$(1)\Rightarrow (3)$. Let $f(z)\equiv 1 \in H^{\infty}$, then $\mathcal {H}(1)(z)=\sum_{n=0}^{\infty}\mu_{n}z^{n}\in Y_{p}$. If $1<q<\infty$, then the desired result follows from Lemma \ref{lem3.2}. If $q=1$, then (\ref{1}) shows that $D^{1}_{0}\subset H^{1}\subset HL(1)$ . This means that $Y_{1}\subset HL(1)$, so we have that $(n+1)^{-1}\mu_{n}\in \ell^{1}$.

$(3)\Rightarrow (2)$. Let $\{f_{k}\}_{k=1}^{\infty}$   be a bounded sequence in  $H^{\infty}$  which converges to $0$ uniformly on every compact subset of  $\mathbb{D}$. Without loss of generality, we may assume   that $f_{k}(0)=0$ for all $k\geq1$ and $\sup_{k\geq 1}||f||_{\infty}\leq 1$.

{\bf Case\textbf{ $1\leq q\leq 2$.}} Since $D^{q}_{q-1} \subset H^{q} \subset HL(q)$, if suffices to prove that $$\lim_{k\rightarrow \infty}||\mathcal {H}_{\mu}(f_{k})||_{D_{q-1}^{q}}=0.$$
 Assume that  $\sum_{n=1}^{\infty}(n+1)^{q-2}\mu_{n}^{q}<\infty$.  Then,
 \[ \begin{split}
\sum_{n=1}^{\infty}(n+1)^{q-2}\mu_{n}^{q}&= \sum_{n=1}^{\infty}\left(\sum_{k=2^{n-1}}^{2^{n}-1}(k+1)^{q-2}\mu_{k}^{q})\right)\\
& \asymp\sum_{n=1}^{\infty}2^{n(q-1)}\mu_{2^{n}}^{q}\\
& \asymp \sum_{n=1}^{\infty}2^{-nq}\left(\sum_{k=2^{n}}^{2^{n+1}-1}(k+1)^{1-\frac{1}{q}}\mu_{k}\right)^{q}.
  \end{split} \]
It follows that
$$ \sum_{n=1}^{\infty}2^{-nq}\left(\sum_{k=2^{n}}^{2^{n+1}-1}(k+1)^{1-\frac{1}{q}}\mu_{k}\right)^{q}<\infty.$$
By Lemma \ref{lem3.1} we have that
 \[ \begin{split}
& \ \ \ \  \int_{0}^{1}(1-r)^{q-1}\left(\sn (n+1)^{1-\frac{1}{q}}\mu_{n}r^{n}\right)^{q}dr \\
  &\asymp \sum_{n=0}^{\infty}2^{-nq}\left(\sum_{k=2^{n}}^{2^{n+1}-1}(k+1)^{1-\frac{1}{q}}\mu_{k}\right)^{q}<\infty.
   \end{split} \]
Therefore, for any  $\varepsilon>0$ there exists a  $0<r_{0}<1$ such that
\begin{equation}\label{6}
\int_{r_{0}}^{1}(1-r)^{q-1}\left(\sn (n+1)^{1-\frac{1}{q}}\mu_{k}r^{n}\right)^{q}dr<\varepsilon.
\end{equation}
It is clear that
\[ \begin{split}
||\mathcal {H}_{\mu}(f_{k})||^{q}_{D_{q-1}^{q}}&=\int_{|z|\leq r_{0}}|\hu(f_{k})'(z)|^{q}(1-|z|)^{q-1}dA(z)+ \int_{r_{0}<|z|<1}|\hu(f_{k})'(z)|^{q}(1-|z|)^{q-1}dA(z)\\
& := J_{1,k}+J_{2,k}.
  \end{split} \]
By the integral representation of  $\hu$, we  get
  \begin{equation}\label{7}
  \hu(f_{k})'(z)=\int_{0}^{1}\frac{tf_{k}(t)}{(1-tz)^{2}}d\mu(t).\end{equation}
Since   $\{f_{k}\}_{k=1}^{\infty}$ is  converge to $0$ uniformly on every compact subset of  $\mathbb{D}$,  for  $|z|\leq r_{0}$ we have that
  \[ \begin{split}
|\hu(f_{k})'(z)|& \leq \int_{0}^{1} \frac{|f_{k}(t)|}{|1-tz|^{2}}d\mu(t)\\
& \lesssim  \int_{0}^{1} |f_{k}(t)|d\mu(t).
    \end{split} \]
It follows Lemma \ref{lem3.3} that
$$J_{1,k}  \rightarrow 0, \  \ \ \mbox{as}\ \ k\rightarrow \infty.$$

 By   Minkowski's inequity and Lemma \ref{lem2.2}, we have that
\[ \begin{split}
M_{q}(r,\hu(f_{k})') &= \left\{\int_{0}^{2\pi}\left|\int_{0}^{1}\frac{tf_{k}(t)}{(1-tre^{i\theta})^{2}}d\mu(t)\right|^{q}d\theta\right\}^{\frac{1}{q}}\\
&\lesssim  \left\{\int_{0}^{2\pi}\left(\int_{0}^{1}\frac{1}{|1-tre^{i\theta}|^{2}}d\mu(t)\right)^{q}d\theta\right\}^{\frac{1}{q}}\\
&\lesssim  \int_{0}^{1}\left(\int_{0}^{2\pi}\frac{d\theta}{|1-tre^{i\theta}|^{2q}}\right)^{\frac{1}{q}}d\mu(t)\\
\end{split} \]
\[ \begin{split}
 &\lesssim \int_{0}^{1}\frac{1}{(1-tr)^{2-\frac{1}{q}}}d\mu(t)\\
&\asymp \sn (n+1)^{1-\frac{1}{q}}\mu_{n}r^{n}.
\end{split} \]

Thus, by the polar coordinate formula and (\ref{6}), we obtain
\begin{align*}
    J_{2,k}& =\int_{r_{0}<|z|<1}|\hu(f_{k})'(z)|^{q}(1-|z|)^{q-1}dA(z)\\
  &   \lesssim \int_{r_{0}}^{1}(1-r)^{q-1}M^{q}_{q}(r,\hu(f_{k})')dr\\
  & \lesssim   \int_{r_{0}}^{1}(1-r)^{q-1}\left(\sn (n+1)^{1-\frac{1}{q}}\mu_{n}r^{n}\right)^{q}dr\\
  & \lesssim \varepsilon.
\end{align*}
Consequently,
$$\lim_{k\rightarrow \infty}||\hu(f_{k})||_{D_{q-1}^{q}}=0.$$

{\bf Case\textbf{ $p>2$.}} \ \  By (\ref{2})  we see that $HL(q)\subset Y_{q}$. To complete the proof, we have to prove that
$\lim_{k\rightarrow \infty}||\hu(f_{k})||_{HL(q)}=0.$

It is clear that the integral $\int_{0}^{1}t^{n}f_{k}(t)d\mu(t)$   converges absolutely for all $n,k \in \mathbb{N}$. It follows that
$$\hu(f_{k})(z)=\int_{0}^{1}\frac{f_{k}(t)}{1-tz}d\mu(t)
= \int_{0}^{1}\sum_{n=0}^{\infty}t^{n}f_{k}(t)z^{n}d\mu(t)
= \sum_{n=0}^{\infty}\left( \int_{0}^{1}t^{n}f_{k}(t)d\mu(t)\right)z^{n}.$$
Since  $\sum_{n=1}^{\infty}(n+1)^{q-2}\mu_{n}^{q}<\infty$, we see that for any $\varepsilon>0$, there exists a positive integer $N$ such that
\begin{equation}\label{8}
\sum_{n=N+1}^{\infty}(n+1)^{q-2}\mu_{n}^{q}<\varepsilon.\end{equation}
For each $k\geq 1$, we have
\[ \begin{split}
&\ \ \ \ \  \sum_{n=0}^{N}(n+1)^{q-2}\left| \int_{0}^{1}t^{n}f_{k}(t)d\mu(t)\right|^{q}\\
&\leq \sum_{n=0}^{N}(n+1)^{q-2}\left(\int_{0}^{1}|f_{k}(t)|d\mu(t)\right)^{q}\\
&\lesssim \left(\int_{0}^{1}|f_{k}(t)|d\mu(t)\right)^{q}.
\end{split} \]
By Lemma \ref{lem3.3}, there there exists $k_{0}\in \mathbb{N}$ such that
\begin{equation}\label{9}
\left(\int_{0}^{1}|f_{k}(t)|d\mu(t)\right)^{q}<\varepsilon \ \ \mbox{for all} \ k>k_{0}.
\end{equation}
Hence, for $k>k_{0}$, by (\ref{8}) and (\ref{9}) we have that
\[ \begin{split}
||\hu(f_{k})||^{q}_{HL(q)}&=\left(\sum_{n=0}^{N}+\sum_{n=N+1}^{\infty}\right)(n+1)^{q-2}\left| \int_{0}^{1}t^{n}f_{k}(t)d\mu(t)\right|^{q}\\
& \lesssim \left(\int_{0}^{1}|f_{k}(t)|d\mu(t)\right)^{q}+ \sup_{k\geq 1}||f_{k}||_{\infty}\sum_{n=N+1}^{\infty}(n+1)^{q-2}\mu_{n}^{q}\\
&\lesssim \varepsilon +\sum_{n=N+1}^{\infty}(n+1)^{q-2}\mu_{n}^{q}\\
& \lesssim\varepsilon.
\end{split} \]
Therefore,
$$\lim_{k\rightarrow \infty} ||\hu(f_{k})||_{HL(q)}=0.$$
The proof is complete.

\end{proof of Theorem 1.4}

\begin{proof of Theorem 1.5}
Let $\{f_{k}\}_{k=1}^{\infty}$   be a bounded sequence in  $H^{\infty}$  which converges to $0$ uniformly on every compact subset of  $\mathbb{D}$.

Since the measure $\mu$ is finite and the integral $\int_{0}^{1}(1-r)^{\frac{1}{q}-2}\log\frac{e}{1-r}dr$ converges, and hence  for any $\varepsilon>0$ there exists a $0<\delta<1$ such that
 $$\int^{1}_{\delta}d\mu(t)<\varepsilon \ \mbox{and}\  \int^{1}_{\delta}(1-r)^{\frac{1}{q}-2}\log\frac{e}{1-r}dr<\varepsilon.$$
For above $\varepsilon>0$, it is clear that  there exists $k_{0}\in \mathbb{N}$ such that
$$\sup_{|w|\leq \delta}|f_{k}(w)|<\varepsilon \ \mbox{for all}\ k>k_{0}.$$
For  $0<r<1$ and $ k>k_{0}$, by Fubini's theorem and Lemma \ref{lem2.2}  we have that
\[ \begin{split}
M_{1}(r,\hu(f_{k}))& =\frac{1}{2\pi}\int_{0}^{2\pi}\left|\int_{0}^{1}\frac{f_{k}(t)}{(1-tre^{i\theta})}d\mu(t)\right|d\theta\\
& \leq  \sup_{|w|\leq \delta}|f_{k}(w)|\int_{0}^{\delta}\frac{1}{2\pi}\int_{0}^{2\pi}\frac{1}{|1-tre^{i\theta}|}d\theta d\mu(t)\\
&+\sup_{k\geq 1}||f_{k}||_{\infty}\int_{\delta}^{1}\frac{1}{2\pi}\int_{0}^{2\pi}\frac{1}{|1-tre^{i\theta}|}d\theta d\mu(t)\\
& \lesssim  \varepsilon \int_{0}^{\delta} \log\frac{e}{1-tr}d\mu(t)+ \int_{\delta}^{1} \log\frac{e}{1-r}d\mu(t)\\
& \lesssim \varepsilon \log\frac{e}{1-r}.
\end{split} \]
Hence, for $k>k_{0}$, we have that
\[ \begin{split}
||\hu (f_{k})||_{B_{q}}&=\int_{0}^{1}(1-r)^{\frac{1}{q}-2}M_{1}(r,\hu(f_{k}))dr \\
& \lesssim \varepsilon \int_{0}^{1} (1-r)^{\frac{1}{q}-2} \log\frac{e}{1-r}dr\\
& \lesssim \varepsilon.
\end{split} \]
This implies that $\hu$ is compact from $H^{\infty}$ to $B_{q}$.\end{proof of Theorem 1.5}

\comment{
Let $f\in H^{\infty}$, then Fubini's theorem and Lemma \ref{lem2.2} show that
\[ \begin{split}
M_{1}(r,\hu(f))& =\frac{1}{2\pi}\int_{0}^{2\pi}\left|\int_{0}^{1}\frac{f(t)}{(1-tre^{i\theta})}d\mu(t)\right|d\theta\\
& \leq ||f||_{\infty} \int_{0}^{\delta}\frac{1}{2\pi}\int_{0}^{2\pi}\frac{1}{|1-tre^{i\theta}|}d\theta d\mu(t)\\
& \lesssim ||f||_{\infty} \int_{0}^{1} \log\frac{e}{1-tr}d\mu(t).
\end{split} \]
Then Fubini's theorem implies that
\[ \begin{split}
\int_{0}^{1}(1-r)^{\frac{1}{q}-2}M_{1}(r,\hu(f))dr &\lesssim  ||f||_{\infty} \int_{0}^{1} (1-r)^{\frac{1}{q}-2} \int_{0}^{1} \log\frac{e}{1-tr}d\mu(t)dr\\
& \leq ||f||_{\infty}\int_{0}^{1} \int_{0}^{1}(1-r)^{\frac{1}{q}-2} \log\frac{e}{1-r}dr d\mu(t)\\
& \lesssim  ||f||_{\infty}.
\end{split} \]
}



\section*{Conflicts of Interest}
The authors declare that there is no conflict of interest.

\section*{Funding}

 The author was supported by  the Natural Science Foundation of Hunan Province (No. 2022JJ30369).



\section*{Availability of data and materials}
Data sharing not applicable to this article as no datasets were generated or analysed during
the current study: the article describes entirely theoretical research.




\pdfbookmark[1]{4 References}{4}

 \end{document}